\documentstyle[12pt]{amsart}
\title[real Lagrangian surfaces with complex tangents]{Divergence of the normalization  for real Lagrangian surfaces near complex tangents}
\author{Xianghong  Gong}
\address{School of Mathematics, Institute for Advanced Study, Princeton, NJ 08540}
\email{ gong@@math.ias.edu}
\newtheorem{thm}{Theorem}[section]

\newtheorem{lemma}[thm]{Lemma}

\begin{document}
%\maketitle
\begin{abstract}
We study real Lagrangian analytic surfaces in ${\Bbb C}^2$ with a 
non-degenerate complex tangent. Webster proved that all such surfaces can be transformed into a quadratic surface by formal symplectic transformations of ${\Bbb C}^2$.
We show that there is a certain dense set of real Lagrangian surfaces which cannot be transformed into the quadratic surface by any holomorphic (convergent) transformation of ${\Bbb C}^2$. The divergence is contributed by the parabolic character of a pair of involutions generated by the real Lagrangian surfaces.
\end{abstract}
\maketitle
\section{Introduction}\setcounter{equation}{0}

We consider a real analytic surface $M$ in ${\Bbb C}^2$. Let $\omega=dz\wedge dp$ be the holomorphic symplectic 2-form on ${\Bbb C}^2$. $M$ is  a {\it real Lagrangian\/} surface if 
\begin{equation}\label{1}
\mbox{Re}\, \omega|_M=0.
\end{equation}
The invariant theory of real Lagrangian submanifolds was initially studied by S.\ ~M.\ ~ Webster~\cite{webster}. It was known that all totally real and real
Lagrangian analytic submanifolds are equivalent under holomorphic symplectic 
transformations.  When $M$ has a non-degenerate complex tangent, Webster proved that under formal symplectic transformations, 
$M$ can be transformed into the quadratic surface
\begin{equation}\label{2}
Q: p=2z\overline z+\overline z^2.
\end{equation}
 Furthermore,  $M$ can be transformed into $Q$
 by holomorphic transformations of ${\Bbb C}^2$ if and only if they are 
equivalent 
through holomorphic symplectic transformations~\cite{webster}.
The purpose of this paper is to  show that there exist real Lagrangian surfaces such that the above normal form  cannot be realized by any holomorphic (convergent)  transformation.

In~\cite{M-W}, J.\ ~K.\ ~Moser and S.\ ~M.\ ~Webster systematically investigated  the holomorphic invariant theory of real surfaces in ${\Bbb C}^2$, where a 
pair of involutions intrinsically attached to the complex tangents plays an important role. We shall see that the divergence for the normalization of the real Lagrangian surfaces is contributed by the parabolic character of the pair of involutions. In~\cite{g2}, a parabolic pair of involutions was also used to show the divergence for the normalization of real analytic glancing hypersurfaces. The main idea of the divergence proof  is inspired by a remark of Moser that divergence of solutions to linearized equations should indicate the same behavior of solutions to the original non-linear equations. In section 3,   we shall derive   a relation between the linearized equations and the original non-linear  ones, which says that  for a certain type of non-linear equations, the existence of the convergent solutions to the linearized equations is indeed a necessary condition for the existence of convergent solutions to the 
non-linear equations.
Our approach comes  directly from the method used for the small divisors (see~\cite{bruno},~\cite{herman} for the references).  In particular, we  follows 
the ideas of  H.\ ~Dulac~\cite{dulac} and  C.\ ~L.\ ~Siegel~\cite{siegel1} 
 closely. 

To state our result, we let $X$ be the set of convergent power series
\begin{equation}\label{3}
r(z,\overline z)=\sum_{i+j>3}r_{ij}z^i\overline z^j,\quad r_{ji}=\overline r_{ij}
\end{equation}
with
\begin{equation}
r_{z\overline z}(z,-z)=0.\label{3+}
\end{equation}
Let us introduce a metric $d$ on  $X$ by
$$
d(r,s)=\sup\left\{|r_{ij}-s_{ij}|^{\frac{1}{i+j}},i+j>3\right\},\quad r,s\in X.
$$
To each $r\in X$, we associate a real Lagrangian surface
\begin{equation}\label{5}
M_r\colon p=2z\overline z+\overline z^2+r_z(z,\overline z).
\end{equation}
Denote by $S$ the set of $r\in X$ such that the corresponding surface $M_r$ cannot be transformed into (\ref{2}) through any holomorphic transformation. We have 
\begin{thm}\label{t1}
$S$ is dense in the metric space $\{ X,d\}$.
\end{thm}
%We mention that in~\cite{webster}, Webster proved that a real Lagrangian 
%surface can be transformed into the quadric $Q$ by holomorphic symplectic 
%transformations if and only if it can be transformed into $Q$ by holomorphic
% transformations.
In section 4, we shall prove  Theorem~\ref{t1}. Using the relation
between the non-linear equations and the linearized equations established in 
section~3, we shall prove that the parabolic pair of involutions generated by 
a real Lagrangian surface is generally not linearizable by any convergent
transformation.

{\bf Acknowledgement}. The main result in this paper is a part of author's thesis~\cite{gt}. The author is grateful to Professor Sidney M.\ ~Webster for the guidance and encouragement. 
The author  acknowledge the support by the NSF grant DMS-9304580 through
 a membership at the Institute for Advanced Study.

\section{Formal theory and linearized equations}\setcounter{equation}{0}\label{section:5.1}
In this section, we shall recall from~\cite{M-W} a pair of   involutions which are intrinsically attached to surfaces with a complex tangent. We also need a result in~\cite{webster} on the formal normalization of a parabolic pair of involutions.

Given a real analytic function $R(z,\overline z)$ with $R(0)=dR(0)=0$, we consider the real Lagrangian surface
$$
M\subset {\Bbb C}^2\colon p=R_z(z,\overline z).
$$
We have
$$
\omega|_M=-R_{z\overline z}dz\wedge d\overline z.
$$
$M$ is totally real if and only if the Levi-form $R_{z\overline z}\neq 0$. It is known that all totally real and real Lagrangian analytic surfaces are equivalent under holomorphic symplectic transformations~\cite{webster}. We now assume that $M$ has a non-degenerate  complex tangent at $0$, i.\ ~e.
\begin{equation}\label{non}
 R_{z\overline z}(0)=0,\quad dR_{z\overline z}(0)\neq 0.
\end{equation}
Then with a suitable change of symplectic coordinates, we may assume that $M$ is given by (\ref{5}) for some  real function (\ref{3}).
From (\ref{non}), we see that $M$ has complex tangents along the smooth curve $C\subset M \colon R_{z\overline z}=0$.  In fact, (\ref{3+}) implies that
$$
C\colon z+\overline z=0.
$$
Here the complex tangents are  {\it parabolic \/} according to E.\ ~Bishop~\cite{bishop}.

Following~\cite{M-W}, we consider the complexification of $M$ defined by
$$
M^c\subset {\Bbb C}^4:\left \{ \begin{array}{l}
p=2zw+w^2+r_z(z,w),\vspace{.5ex}\\
q=2zw+z^2+\overline r_{\overline z}(w,z).
\end{array} \right.
$$
We shall use $(z,w)$ as the coordinates to identify $M^c$ with ${\Bbb C}^2$. Consider the projection $\pi_1\colon (z,p,w,q)\to (w,q)$.  The restriction of $\pi_1$ to $M^c$ is a 
double-sheeted branched covering, and it induces a covering transformation $\tau_1\colon M^c\to M^c$.  Notice that $w$ and $q$ are invariant under $\tau_1$. Then $\tau_1\colon (z, w)\to(z^\prime,w^\prime)$ is implicitly defined by
\begin{equation}                    \label{eq:5.1.3}
\tau_1:\left\{ \begin{array}{l}
z^{\prime}=-z-2w-{1\over z^\prime-z}\{ r_{\overline z}(z^\prime,w)-r_{\overline z}(z,w)\},
\vspace{.5ex}\\ w^\prime=w.  
\end{array}\right. \end{equation}
We also have the projection $\pi_2\colon (z,p,w,q)\to(z,p)$ inducing another involution $\tau_2$ of $M^c$.  We notice that $M^c$ is invariant under
the complex conjugation  
$$ (z,p,w,q)\to (\overline w,\overline q,\overline z,\overline p).$$
Hence, its restriction to $M^c$ is an anti-holomophic involution $\rho\colon (z,w)\to(\overline w,\overline z)$. From the relation $\pi_2=c\pi_1\rho$ for $c(
z,p)=(\overline z,\overline p)$, it follows that $\tau_1$ and $\tau_2$ satisfy the reality condition
\begin{equation}            \label{eq:5.1.8}
\tau_2=\rho\circ\tau_1\circ\rho.
\end{equation}

Introduce the following coordinates for $M^c$
\begin{equation}
x=z+w,\hspace{2em} y=z-w.           \label{eq:5.1.4}
\end{equation}
Then $\rho$ takes the form
\begin{equation}                    \label{eq:5.1.5}
\rho (x,y)=(\overline x,-\overline y).
\end{equation}
Now the pair of involutions $\tau_1, \tau_2$ can be written as
\begin{equation}
            \label{eq:5.1.6}
\tau_j(x,y)=\tau_j^*(x,y)+H_j(x,y),\quad  H_j(x,y)=O(2)
\end{equation}
with
\begin{equation}\label{7}
\tau_1^*(x,y)=(-x,-2x+y), \quad \tau_2^*(x,y)=(-x,2x+y).
\end{equation}
 Since $\tau_j^2=id$, then
\begin{equation}\label{7+}
\tau_j^*\circ H_j+H_j\circ\tau_j= 0 .
\end{equation}
We also notice that the branch points of $\pi_j$, i.\ ~e.\ the fixed points of $\tau_j$, are given by
$$2z+2w+H_{z\overline z}(z,w)=0,$$
i.\ ~e,\ $x=0$. Hence, we have
\begin{equation}\label{7++}
\quad \ H_j(0,y)= 0.
\end{equation}
  The reality condition on $\{\tau_1, \tau_2\}$ is still given by (\ref{eq:5.1.8}).

%Let $\omega$ be the restriction of the 2-form $dz\wedge dp$ to $M^c$. Then 
% it is known that $\{ \tau_j, \rho,\omega\}$ determines intrinsically the real% Lagrangian surfaces~\cite{webster}. More specifically, if
From~\cite{M-W}, we now recall an intrinsic property of the pair of involutions generated by a surface with a non-degenerate complex tangent as follows. Let $\{ \hat\tau_j, \rho \}$ and $\{ \tau_j, \rho \}$ be two pairs of involutions corresponding to two real (Lagrangian) analytic surfaces $\hat M$ and $M$ in the form (\ref{5}). Then
  $M$ and $\hat M$ are equivalent through holomorphic transformations if and only if there is a biholomorphic mapping $\Phi\colon M^c\to \hat M^c$ such that
$$
\Phi^{-1}\hat{\tau}_j	\Phi=\tau_j,\quad %\ \Phi^*\hat\omega=\omega,\
 \rho\Phi=\Phi\rho.
$$
Notice that involutions generated by the quadratic surface (\ref{2}) are the linear involutions  (\ref{7}). We shall prove that there is a dense subset $S$ in $X$ of which the involutions generated by the corresponding surfaces are not linearizable by any convergent transformations of $M^c$, from which Theorem~\ref{t1} follows immediately.

 A convergent or formal transformation
\begin{equation}
              \label{eq:5.1.9}
   \Phi\colon x\mapsto x+u(x, y),\quad y\mapsto y+v(x, y)
\end{equation}
is said to be {\it normalized\/} if
\begin{equation}\label{9+}
  u(0,y) = 0,\quad
         u(x,0)=u(-x,0),\quad  v(x,0)=- v(-x,0).
\end{equation}

\begin{lemma}[\cite{webster}]\label{l1}
 Let $\tau_1$ and $\tau_2$ be a pair of involutions defined by $(\ref{eq:5.1.6}), (\ref{7})$ and $(\ref{7+}).$ Then there exists a unique normalized formal transformation $\Phi$ such that $\Phi\tau_j\Phi^{-1}=\tau_j^*$ for $j=1,2.$
\end{lemma}
We also need the following.
\begin{lemma}[\cite{g2}]\label{l2}
Let $\tau_j\ (j=1,2)$ be as in Lemma~\ref{l1}$.$ Then $\tau_j$ are linearizable by convergent transformations if and only if the above unique normalized transformation is convergent$.$
\end{lemma}

%\subsection{Linearized equations}\setcounter{equation}{0}
 We now want to discuss the linearized equations for the pair of involutions. Put
  $$
H_j=( f_{j},  g_{j}),\quad  j=1,2 .
$$
 Consider the composition
\begin{equation} 
\sigma=\tau_1\tau_2: \left\{           \label{eq:5.1.10}
\begin{array}{l}
x^\prime  =  x +G_1(x,y),\vspace{.5ex}\\
y^\prime  =  4x+y +G_2(x,y),
\end{array}\right.
\end{equation}
where
\begin{equation}\label{g1}
G_1=
- f_2 + f_1\circ\tau_2 ,\quad G_2=
-2 f_2  + g_2  + g_1\circ\tau_2.
\end{equation}
From $\Phi\tau_j=\tau_j^*\Phi\ (j=1,2)$, we get
\begin{equation}
\begin{array}{l}                         \label{eq:5.1.11}
 u\circ\sigma- u=-G_1,\vspace{.5ex}\\
 v\circ\sigma- v=-G_2+4 u.
\end{array}
\end{equation}          
This leads us to the following linearized equation
\begin{equation} \label{eq:5.1.12}
 u(x, 4x+y)-u(x,y)=-G_1(x,y).
\end{equation}
From (\ref{7++}), it follows that $x|G_1(x,y)$. Set
\begin{equation}                          \label{eq:5.1.13}
\begin{array}{l}
\displaystyle{ K=\sum_{k=0}^\infty 4^k\beta_kx^kD^k,\quad D=  {\partial\over \partial y}}, \vspace{.5ex}\\
E(z)=\displaystyle{ {z\over e^z-1}=\sum_{k=0}^\infty \beta_kz^k}.
\end{array}
\end{equation}
Rewrite (\ref{eq:5.1.11}) as
$$ (e^{4xD}-1) u(x,y)=-G_1(x,y). $$
Clearly, $K(e^{4xD}-1)=4xD$. By applying $K$ to the above, we finally reduce (\ref{eq:5.1.12}) to
\begin{equation}                        \label{eq:5.1.14}
\partial_y  u(x,y)={1\over 4}K\, a(x,y),\quad
G_1(x,y)=-xa(x,y).
\end{equation}
From (\ref{eq:5.1.13}), it  follows  that the linearized equation   (\ref{eq:5.1.12}) has only divergent solutions. This can be seen easily if   $G_1$ can be  arbitrarily chosen. However,  as $G_1$ comes from real Lagrangian surfaces, we  need to study the linearized equations more closely.

\section{ Linearizations}\setcounter{equation}{0}
\label{section:5.2}

In this section, we shall investigate the relation between non-linear equations and their first order approximate equations, i.\ ~e.\ ~ the linearized equations. We shall consider 
non-linear equations which are formally solvable. Then under suitable conditions, we shall formally solve the linearized equations,  and show that the existence of a divergent solution to linearized equations implies that the original non-linear equations have neither convergent solutions.

Consider a system of equations
\begin{equation}\label{f}
F(x,y)=0, \qquad (x,y)\in X\times Y,
\end{equation}
where $X=\oplus_{k=1}^\infty X_k, \ Y=\oplus_{k=1}^\infty Y_k$, and
$$
F\colon X\times Y\to Z,\quad F(0,0)=0
$$
with $Z=\oplus_{k=1}^\infty Z_k$. We shall denote by $\pi_k$ the projection $X\to X_k$, as well as the projections $Y\to Y_k$ and $Z\to Z_k$.  For each $x\in X$, there is a formal decomposition $x=x_1+x_2+\ldots$ with $x_k\in X_k$. We 
put
$$
[x]_k=x_1+x_2+\ldots+ x_k
$$
with $[x]_0=0$.
 Assume that $X_k,Y_k,Z_k$ are finitely dimensional real vector spaces identified with Euclidean spaces. We then introduce the product topology on $X,Y,Z$. Thus, a sequence $\{x^{(n)}\}\subset X$ is convergent if and only if $\{x_k^{(n)}\}_{n=1}^\infty$ converges in $X_k$ for all $k$.

The Fr\'echet derivative of $F$ at $(0,0)$ is defined by
$$
DF(x,y)=\lim_{{\Bbb R}\ni t\to 0}\frac{F(tx,ty)-F(0,0)}{t},\quad x\in X,\  y\in Y. 
$$
We assume that $DF\colon X\times Y\to Z$ is a well-defined linear mapping. We further assume that $DF$ is {\it homogeneous}, i.\ ~e.
\begin{equation}\label{h}
DF\colon X_k\times Y_k\to Z_k.
\end{equation}
Decompose
$$
DF(x,y)=D_1F(x)+D_2F(y)
$$
with $D_1F(x)=DF(x,0) $ and $D_2F(y)=DF(0,y)$. Put
$$
QF=F-DF.
$$
 We assume that for $k\geq 1$
\begin{equation}\label{q}
\pi_kQF(x,y)=\pi_kQ([x]_{k-1},[y]_{k-1}).
\end{equation}

With the above notations and assumptions, we have
\begin{lemma}\label{l3} Let $F$ be as above$.$ Assume that
 $F$ 
satisfies the following conditions$:$
\newcounter{l1}
\begin{list}
{$(\roman{l1})$}{\usecounter{l1}\setlength{\rightmargin}{\leftmargin}}
\item
There is a solution operator $P\colon X\to Y$ with $F(x,P(x))=0$ and $P(0)=0.$
\item
$D_2F\colon Y\to Z$ is injective$.$
\end{list}
Then  $DF(x,y)=0$ is solvable by $L\colon X\to Y$ with
\begin{equation}\label{dl}
L(x)=\sum_{k=1}^\infty \pi_kP(x_k),\quad x\in X.
\end{equation}
Moreover$,$ $\widetilde P=P-L$ satisfies that
\begin{equation}\label{tq}
\pi_k\widetilde P(x)=\pi_k\widetilde P([x]_{k-1}).
\end{equation}
\end{lemma}
\pf Fix $x\in X$.  Notice that $L$ is homogeneous. Hence, it suffices to show that
(\ref{tq}) holds and
\begin{equation}\label{*}
DF(x_k, L(x_k))=0.
\end{equation}
From (\ref{q}), we see immediately that (\ref{*}) holds for $k=1$. We also have
$$
DF(x_1, \pi_1P(x))=0.$$
Since $D_2F$ is injective, then we get $\pi_1\widetilde P(x)=0$. We now assume that  (\ref{tq}) and (\ref{*}) 
 hold for $k<n$. From (\ref{tq}), we get
$[P(x)]_{n-1}=[P([x]_{n-1})]_{n-1}$. In particular, $[P(x_n)]_{n-1}=0$.
Now (\ref{q}) gives
$$
DF(x_n, \pi_n P(x_n))=0.
$$
We also have
$$
DF(x_n,\pi_nP(x))=-\pi_nQF(x,P(x)).
$$
Since $[P(x)]_{n-1}$ depends only on $[x]_{n-1}$, we obtain from the last two identities that
$$
DF(0,\pi_n\widetilde P(x))=-\pi_nQF\left([x]_{n-1},[P([x]_{n-1})]_{n-1}\right).
$$
Since $ D_2F$ is injective, we can find a left inverse $K\colon Z\to Y_0$ of $D_2F$ such that $K\colon Z_k\to  Y_k$. Thus, we get
$$\pi_n\widetilde P(x)=-K\pi_nQF\left([x]_{n-1},[P([x]_{n-1})]_{n-1}\right).
$$
This proves that $\pi_n\widetilde P(x)$ depends only on $[x]_{n-1}$. Therefore, (\ref{tq}) holds for $k=n$. The proof of Lemma~\ref{l3} is complete.

Assume that $X$ is endowed with a metric $d$ such that
$$d(x^\prime,x^{\prime\prime})=\sup\{d(x_k^\prime,x_k^{\prime\prime}); k=1,2,\ldots\}.
$$
We also assume that each $Y_k$ is endowed with  a metric $\tilde d_k$ which is invariant under translation. Put
$$
\tilde d(y^\prime,y^{\prime\prime})=\sup\{\tilde d_k(y_k^\prime,y_k^{\prime\prime}); k=1,2,\ldots\}, \quad y^\prime, y^{\prime\prime}\in Y,
$$
and $\hat Y=\{ y\in Y;\tau(y,0)<\infty\}$.

 With the above notations and assumptions, we want to prove the following.
\begin{lemma}\label{l4}
Suppose that there exits $x^*\in X$ with $LP(x^*)\notin \hat Y.$ Let $\epsilon_0=d(x^*,0).$ Then for any $x\in X,$ there is $x^\prime\in X$ with $d(x,x^\prime)\leq\epsilon_0$ such that $P(x^\prime)\notin\hat Y.$
\end{lemma}

For the application of Lemma~\ref{l4} to real Lagrangian surfaces, $Y$ will be taken as a certain space of normalized formal solutions, and the metric on $Y$ will be so chosen that $\hat Y$ is precisely the set of convergent solutions among all the normalized formal solutions. Thus, Lemma~\ref{l4} says that the existence of convergent solutions to linearized equations is a necessary condition for the existence of a convergent solution to the original functional equations. On the other hand, we notice that Siegel's  theory on the Hamiltonian systems~\cite{siegel2} concluded that there are functional equations of which the convergent solutions  does exist for the linearized equations, but not for the original functional equations. Finally, we mention that
the metric which we put on the space of convergent power series in Theorem~\ref{t1} is weaker  than that used by Siegel in~\cite{siegel1}, where the small divisors are essential.

{\it Proof of \mbox{Lemma~\ref{l4}}}. Since $LP(x^*)\notin \hat Y$, then we can choose a sequence of positive integers $n_k$ such that 
\begin{equation}\label{le}
\tilde d_{n_k}\left(LP(x_{n_k}^*),0\right)\geq k
\end{equation}
for all $k$. Fixing $x\in X$, we put $x_n^\prime=x_n$ for all $n\neq n_k$. To determine $\{x_{n_k}^\prime\}$, we assume that for $m<n_k$, all $x_m$ have been 
so chosen  that
\begin{equation}\label{j/2}
\tilde d_{n_j}(\pi_{n_j}P([x^\prime]_{n_j}),0)\geq j/2
\end{equation}
for all $j<k$. Let $\tilde x=[x^\prime]_{{n_k}-1}+x_{n_k}$. If 
$$
\tilde d_{n_j}(p_{n_k}(\tilde x),0)\geq j/2,
$$
we put $x_{n_k}^\prime=x_{n_k}$. Then (\ref{j/2}) holds for $j=k$. Otherwise, we choose
$$
x_{n_k}^\prime= x_{n_k}+x_{n_k}^*.
$$
Then from Lemma~\ref{l3}, it follows that
$$
\pi_{n_k}P([x^\prime]_{n_k})=LP(x_{n_k}^*)+\pi_{n_k}P(\tilde x).
$$
Since all $\tilde d_n$ are invariant under translation, then
$$
\tilde d_{n_k}(\pi_{n_k}P([x^\prime]_{n_k},0)\geq \tilde d_{n_k}(Lp(x_{n_k}^*),0)-\tilde d_{n_k}P(\tilde x,0).
$$
Now (\ref{le}) implies that (\ref{j/2}) holds for $j=k$. Thus, we have chosen $x^\prime$ such that (\ref{j/2}) holds for all $j$. Notice that $[P(x)]_k$ depends only on $[x]_k$ for all $k$. Therefore,  (\ref{j/2}) implies that 
$P(x^\prime)\notin\hat Y$. The proof of Lemma~\ref{l4} is complete.

\section{Proof of Theorem~\ref{t1}}\setcounter{equation}{0}

In this section, we shall give a proof of Theorem~\ref{t1} by using Lemma~\ref{l4}. Thus, we shall verify that for the problem of linearizing the pair of involutions generated by a real Lagrangian surface, its linearized equations have divergent solutions..

We first give some notations. Let $(X, d)$ be as in Theorem~\ref{t1}. Denote by $X_k$ the set of homogeneous polynomials $r\in X$ with degree $k+2$. 
 Let $Z_k\ (k\geq 2)$  be the set of ordered
 power series $(w_1,w_2,w_3,w_4)$  in $x,y$, where each $w_j$ is a homogeneous polynomial of degree $k$.
Put $Z$ to be the direct sum of $Z_k\ (k\geq 2)$.
 Denote by $Y_k$ the set of ordered pairs $(u,v)$  of homogeneous polynomials
 of degree $k$ in $x,y$ which satisfy the normalizing condition (\ref{9+}). Let $Y$ be the direct sum of $Y_k\ (k\geq 2)$. The metric on $Y_k$ is defined by 
$$
\tilde d_k\left((u^\prime,v^\prime),(u,v) \right)=
\max_{i+j=k}\left\{
|u_{ij}^\prime-u_{ij}|^{\frac{1}{i+j}},|v_{ij}^\prime
-v_{ij}|^{\frac{1}{i+j}}
\right\}.
$$

We now put the system of equations $\Phi\tau_j\Phi^{-1}=\tau_j^*\ (j=1,2)$ into the form 
 (\ref{f}) with 
$$
F(r,u,v)=(\Phi\circ\tau_1-\tau_1^*\circ\Phi, \Phi\circ\tau_2-\tau_2^*\circ\Phi).
$$
Then we have
$$
D_2F(u,v)=(u\circ\tau_1^*+u,v\circ\tau_1^*-v+2u,u\circ\tau_2^*+u,v\circ\tau_2^*-v-2u).
$$
In order to apply Lemma~\ref{l4}, we need to verify that $D_2F$ is injective. This is essentially contained in the argument in~\cite{webster}. To give the details, 
we notice that $D_2F(u,v)=0$ implies that $u$ is invariant under both $\tau_1^*$ and $\tau_2^*$. Then $u(x,y)$ depends only $x$, and it contains only power of $x$ with odd order (see Lemma 4.1 in~\cite{g2}). From the normalizing condition (\ref{9+}),
it  follows that $u=0$. We now know that $v$ is skew-invariant under $\tau_j$ 
$(j=1,2)$. Then $v(x,y)$ depends only $x$ and contains only power of $x$ with
 even order. Thus, (\ref{9+}) implies that $v=0$.    It is easy to verify that $DF$ is homogeneous and $F$ satisfies (\ref{q}).

To compute $D_1F$, we fix $r\in X$ and put
\begin{equation}                     \label{eq:5.2.2}
\tau_1:\left\{ \begin{array}{l}
z^\prime=-z-2w+q(z,w), \quad q(z,w)=O(2),\vspace{.5ex}\\
w^\prime=w.\end{array}\right. 
\end{equation}
By the implicit formula (\ref{eq:5.1.3}), we obtain
$$
q(z,w)=Lq(z,w)+Qq(z,w)
$$
with
\begin{equation}        \label{eq:5.2.3}
Lq(z,w)=-\frac{1}{2z+2w}\left\{r_{\overline z}(-z-2w,w)-r_{\overline z}(z,w)
\right\},
\end{equation}
where
 each coefficient 
$(Qq)_{i,j}$ of $Qq$ does not contain  linear terms in $r$ and terms $r_{i^\prime,j^\prime}$ with $i^\prime+j^\prime\geq i+j+2$. From (\ref{eq:5.1.4}), (\ref{eq:5.1.6}) and (\ref{eq:5.2.2}), we  get
\begin{equation}
f_1(x,y)=g_1(x,y)              \label{eq:5.2.4}
=q(\frac{x+y}{2},\frac{x-y}{2}).
\end{equation}
Since $\tau_2=\rho\tau_1\rho$, we have
$$\tau_2(x,y)=\left(-x+\overline{f}_1(x,-y),2x+y-\overline{g}_1(x,-y)\right).$$
It follows from (\ref{eq:5.2.4}) that
\begin{equation}
f_2(x,y)=-g_2(x,y)              \label{eq:5.2.5}
=\overline{q}(\frac{x-y}{2},\frac{x+y}{2}).
\end{equation}

We now can obtain
$$
D_1F(r)=\left(Lq(\frac{x+y}{2},\frac{x-y}{2}),Lq(\frac{x+y}{2},\frac{x-y}{2}),\overline{Lq}(\frac{x-y}{2},\frac{x+y}{2}),-\overline{Lq}(\frac{x-y}{2},\frac{x+y}{2})\right).
$$
The equation $DF(r,u,v)=0$ implies that
\begin{gather*}
u\circ\tau_1^*(x,y)+u(x,y)=-Lq(\frac{x+y}{2},\frac{x-y}{2}),\\
u\circ\tau_2^*(x,y)+u(x,y)=-\overline{Lq}(\frac{x-y}{2},\frac{x+y}{2}).
\end{gather*}
This leads to 
\begin{equation}\label{a}
u
(x,4x+y)-
u
(x,y)
=A(x,y)
\end{equation}
with
$$
A(x,y)=\overline{Lq}(\frac{x-y}{2},\frac{x+y}{2})-Lq(\frac{x+y}{2},\frac{3x+y}{2}).
$$
Then from (\ref{eq:5.2.3}), we get
\begin{eqnarray*} \label{eq:5.2.6}
2xA(x,y)&= &r_{\overline z}(\frac{5x+y}{2},-\frac{3x+y}{2})-
                        r_{\overline z}(\frac{x+y}{2},- \frac{3x+y}{2})\\
          & &\quad +\, \overline{r}_{\overline z}(-\frac{3x+y}{2},\frac{x+y}{2})-
                 \overline{ r}_{\overline z}(\frac{x-y}{2},\frac{x+y}{2}).
\end{eqnarray*}

Applying (\ref{eq:5.1.14}), we have  
\begin{equation}       \label{eq:5.2.7}
x\partial_yu
(x,y)=\frac{1}{4}KA(x,y).
\end{equation}
Let
\begin{equation}          \label{eq:5.2.8}
e_{n+2}=\epsilon^{n+2}\left(i^{n-1}z^n\overline z^2+(-i)^{n-1}\overline{z}^n z^2\right) \equiv i^{n-1}P_{n+2}(z,\overline z).
\end{equation}
Then $e_{n+2}\in X_n$ and $\overline{e}_{n+2}(z,\overline z)\equiv -i^{n-1}P_{n+2}(-z,-\overline z)$. Let $A_n$ be given by (\ref{eq:5.2.6}) in which $r$ is replaced by $e_{n+2}$. Then
\begin{equation}
\begin{array}{l}           \label{eq:5.2.9}
2^{n+2}i^{1-n}xA_n(x,y)=P_{n+2,\overline z}(5x+y,-3x-y)\vspace{.5ex}\\
 \qquad -P_{n+2,\overline z}(x+y,-3x-y)+P_{n+2,\overline z}(3x+y,-x-y)\vspace{.5ex}\\ \qquad -P_{n+2,\overline z}(-x+y,-x-y).
\end{array}
\end{equation}
We now put
$$
r^*(z,\overline z)=\sum_{n=4}^\infty e_n(z,\overline z),\quad P(z,\overline z)=\sum_{n=4}^\infty P_n(z,\overline z).
$$
Then 
we have
$$
A(x,y)=\sum_{n=2}^\infty A_n(x,y).$$
\begin{lemma}\label{lemma:5.2.1}
$KA(x,y)$ diverges.
\end{lemma}
\pf We need to compute
$$
K\mid_{y=0}P_{n+2,\overline z}(ax+y,bx-y)\equiv \gamma_{n+1}x^{n+1}.
$$
We have
$$P_{n+2,\overline z}(z,\overline z)=\epsilon^{n+2}\left(2z^n\overline z+(-1)^{n-1}n\overline z^{n-1}z^2\right).$$
This gives
\begin{eqnarray*}
\frac{1}{\epsilon^{n+2}}\gamma_{n+1}&=&\frac{1}{\epsilon^{n+2}}\sum_{j=0}^{n+1}4^j\beta_j(\frac{\partial}{\partial y})^j\mid_{x=1,y=0}P_{n+2,\overline z}(ax+y,bx-y)\\
&=&2\sum_{j=0}^n 4^j\beta_j\frac{n!}{(n-j)!}a^{n-j}b
-2\sum_{j=1}^{n+1} 4^j\beta_j\frac{n!}{(n-j+1)!}a^{n-j+1}\\
& & +(-1)^{n-1}n\sum_{j=0}^{n-1} 4^j\beta_j\frac{(n-1)!}{(n-1-j)!}(-1)^jb^{n-1-j}a^2\\
& &+(-1)^{n-1}n\sum_{j=1}^n 4^j\beta_j2\binom j1 \frac{(n-1)!}{(n-j)!}(-1)^{j-1}b^{n-j}a\\
& &+(-1)^{n-1}n\sum_{j=2}^{n+1} 4^j2\binom j2\beta_j\frac{(n-1)!}{(n+1-j)!}(-1)^{j-2}b^{n+1-j}.
\end{eqnarray*}
Using Cauchy product, one can obtain
\begin{eqnarray*}
\sum_{n=2}^\infty \frac{\gamma_{n+1}}{n!\epsilon^{n+2}}x^{n+1}
&=&2bxE(x)e^{ax}-2E(x)e^{ax} +a^2x^2E(x)e^{-bx}\\
& &+2ax^2E^\prime(x) e^{-bx}+xE^{\prime\prime}(x) e^{-bx}\\
& \equiv & 2ax^2E^\prime(x)
 e^{-bx}
+\tilde{S}_{a,b}(x)\equiv S_{a,b}(x).
\end{eqnarray*}
Thus from (\ref{eq:5.2.9}), we obtain
$$
\sum_{n=2}^\infty\frac{i^{1-n}2^{n+2}}{n!\epsilon^{n+2}}x(KA_n)(x,0)=S(x)
$$
with
$$
S(x)=S_{5,-3}(x)-S_{1,-3}(x)+S_{3,-1}(x)-S_{-1,-1}(x).
$$
 For a non-zero integer $k$,  
$$
\tilde S(x)=\tilde S_{5,-3}(x)-\tilde S_{1,-3}(x)+\tilde S_{3,-1}(x)-\tilde S_{-1,-1}(x)
$$
has no pole of order $2$ at $x=k\pi i/2$. Hence $S(x)$ has a pole of order $2$ at $x=k\pi i/2$ if $e^{3x}+e^x\neq 0$, i.\ e.\ if $k$ is even. Now let
$S(x)=\sum_{n=2}^\infty S_nx^n$. Then we have $\limsup_{n\to \infty}\root n \of{ |S_n|}>0.$ Therefore, one can see that  
$
(KA)(x,0)$
diverges.  The proof of Lemma~\ref{lemma:5.2.1} is complete.

We have proved that (\ref{a}) has a divergent solution. Hence, the equation
 $DF(r^*,u,v)=0$ has a divergent solution for $(u,v)$. Since $d(r^*,0)=\epsilon$, then Lemma~\ref{l4} implies that for any $r\in X$, there is $r^\prime\in X$ with $d(r^\prime,0)\leq \epsilon$ such that the equation $F(r^\prime,u,v)$ has no convergent solution in $Y$. From Lemma~\ref{l3}, it follows that the corresponding pair of involutions $\{\tau_1,\tau_2\}$ is not linearizable by any convergent transformations. This proves Theorem~\ref{t1}.

\bibliographystyle{siam}

\end{document}